\documentclass[showpacs,preprintnumbers,amsmath,amssymb]{revtex4}

\usepackage{graphicx}
\usepackage{bm}
\usepackage{epstopdf}
\usepackage{float}
\usepackage{dcolumn}
\usepackage{bm}
\usepackage{mathrsfs}
\usepackage{amsmath}

\begin{document}


\title{New indefinite integrals of confluent Heun functions}

\author{D. Batic}
\email{davide.batic@ku.ac.ae}
\affiliation{%
Department of Mathematics,\\  Khalifa University of Science and Technology,\\ 
Abu Dhabi, United Arab Emirates 
} 
\author{O. Forrest}
\email{omar_forrest@yahoo.com}
\affiliation{%
Department of Mathematics,\\  University of the West Indies,\\
Kingston 6, Jamaica 
}
\author{M. Nowakowski}
\email{mnowakos@uniandes.edu.co}
\affiliation{
Departamento de Fisica,\\ Universidad de los Andes,\\ 
Cra.1E No.18A-10, Bogota, Colombia
}%

\date{\today}

\begin{abstract}
We obtain some new formulae to compute the first derivative of confluent and biconfluent Heun functions under the minimal assumption of fixing only one parameter. These results together with the Lagrangian formulation of a general homogeneous linear ordinary differential equation allow to construct  several new indefinite integrals for the confluent, biconfluent, doubly confluent, and triconfluent Heun functions.
\end{abstract}

\pacs{XXX}
\maketitle

\section{Introduction}
It is well-known that the ordinary second-order linear differential equation
\begin{equation}\label{0}
\mathcal{L}(y)(x)=0,\quad\mathcal{L}:=\frac{d^2}{dx^2}+p(x)\frac{d}{dx}+q(x),
\end{equation} 
admits the indefinite integral \cite{con1,con2}
\begin{equation}\label{1}
\int\! f(x)\left[h^{''}(x)+p(x)h^{'}(x)+q(x)h(x)\right]y(x)\,\mathrm{d}x=f(x)W(y,h)(x)+c,
\end{equation}
where
\begin{equation}\label{f}
f(x):=e^{\int\! p(x)\,\mathrm{d}x},
\end{equation}
$h$, $p$, and $q$ are complex-valued differentiable functions in the variable $x\in\mathbb{R}$ with $h$ at least twice continuously differentiable, and $W$ denotes the Wronskian. There are quite a few methods at our disposal to get an indefinite integral involving the solution of (\ref{0}) by using (\ref{1}). For instance, we can express $h$ in terms of a simple elementary function such as 
\begin{equation}\label{scelta_h}
h(x)=x^m e^{\rho x^\ell}\left\{\begin{array}{c}
\sin{(k_1 x)}\\
\cos{(k_2 x)}
\end{array}
\right\}
\end{equation}
with $m,\ell\in\mathbb{N}_0:=\mathbb{N}\cup\{0\}$, $\rho,k_1,k_2\in\mathbb{C}$, or when $h$ is a constant function, (\ref{1}) takes the simpler form
\begin{equation}\label{zero}
\int\! f(x)q(x)y(x)\,\mathrm{d}x=-f(x)y^{'}(x)+c.
\end{equation}
If, instead, we look for an integral involving $y(x)$ only, it is desirable to choose $h$ so that
\begin{equation}\label{y_alone}
h^{''}(x)+p(x)h^{'}(x)+q(x)h(x)=\frac{1}{f(x)}.
\end{equation}
Moreover, we may also wish to specify $h$ as a solution to (\ref{0}) with one or two  terms in (\ref{0}) deleted. Typical choices are 
\begin{eqnarray}
h^{''}(x)+p(x)h^{'}(x)&=&0,\label{2}\\
p(x)h^{'}(x)+q(x)h(x)&=&0,\label{3}\\
h^{''}(x)+q(x)h(x)&=&0.\label{4}
\end{eqnarray}
In the case that $p$ and $q$ consists of multiple terms we might also try to specify $h$ with $p$ or $q$ with some of their subterms removed. Another possibility of getting interesting integrals from (\ref{1}) emerges from choosing $h$ to be a particular solution of the equation conjugate to (\ref{0}), that is $h^{''}(x)+p(x)h^{'}(x)+\overline{q}(x)h(x)=0$. This equation has the same $p$ as in (\ref{0}) but a different $q$. If we recall that a transformation of the dependent variable can always be set up to to make any two differential equations conjugate, it is clear that this method may lead to the construction of indefinite integrals containing products of the solutions of the two conjugate equations. More precisely, we have \cite{con1,con2}
\begin{equation}\label{intcon}
\int\! f(x)\left[q(x)-\overline{q}(x)\right]h(x)y(x)\,\mathrm{d}x=f(x)\left[h^{'}(x)y(x)-h(x)y^{'}(x)\right]+c.
\end{equation}
In a companion paper \cite{us}, we constructed new indefinite integral involving Heun functions. Here, we focus on the derivation of indefinite integrals containing confluent, biconfluent, doubly confluent, and triconfluent Heun functions \cite{ronv}. All the aforementioned functions are solutions of the differential equation (\ref{0}) with $p$ and $q$ specified as in Table~\ref{Heun}.
\begin{table}[ht]
\centering
\caption{Schematic representation of the Heun family of differential equations. All parameters otherwise stated are complex. In the case of the Heun operator $\epsilon=\alpha+\beta+1-\delta-\gamma$ and $a\in\mathbb{R}\backslash\{0,1\}$.\cite{ronv,Dec}}
\label{Heun}
\begin{tabular}{lll}
\hline
Equation &  $p(x)$ & $q(x)$ \\ \hline
Heun     &  \begingroup\makeatletter\def\f@size{10}\check@mathfonts$\frac{\gamma}{x}+\frac{\delta}{x-1}+\frac{\epsilon}{x-a}$\endgroup   & \begingroup\makeatletter\def\f@size{10}\check@mathfonts$\frac{\alpha\beta x-q}{x(x-1)(x-a)}$\endgroup\\
confluent Heun       & \begingroup\makeatletter\def\f@size{10}\check@mathfonts$\alpha+\frac{1+\beta}{x}+\frac{1+\gamma}{x-1}$\endgroup   & \begingroup\makeatletter\def\f@size{10}\check@mathfonts $\frac{[(\beta+\gamma+2)\alpha+2\delta]x-(\beta+1)\alpha+(\gamma+1)\beta+2\eta+\gamma}{2x(x-1)}$ \endgroup \\             
biconfluent Heun & \begingroup\makeatletter\def\f@size{10}\check@mathfonts$\frac{\alpha+1}{x}-\beta-2x$\endgroup &\begingroup\makeatletter\def\f@size{10}\check@mathfonts$\gamma-\alpha-2-\frac{\delta+\beta(\alpha+1)}{2x}$\endgroup \\
doubly confluent Heun &\begingroup\makeatletter\def\f@size{10}\check@mathfonts$\frac{2x^3-\alpha x^2-2x-\alpha}{(x-1)^2(x+1)^2}$\endgroup & \begingroup\makeatletter\def\f@size{10}\check@mathfonts$\frac{\beta x^2+(\gamma+2\alpha)x+\delta}{(x-1)^3(x+1)^3}$\endgroup\\
triconfluent Heun &\begingroup\makeatletter\def\f@size{10}\check@mathfonts$-(\gamma+3x^2)$\endgroup &\begingroup\makeatletter\def\f@size{10}\check@mathfonts$\alpha+(\beta-3)x$\endgroup\\ \\ \hline
\end{tabular}
\end{table} 
All findings presented in the subsequent sections have been checked by means of Maple. This is extremely important, as it means our results are immediately applicable, and need not await further analysis giving formulas for the derivatives of the doubly confluent and triconfluent Heun functions. Last but not least, we believe the present work represents a substantial first attack on the problem of integrating the confluent Heun functions considered in this paper.

\section{Examples for some special functions}
Since the function $h$ appearing in (\ref{1}) is not fixed, we have in practice an unlimited number of cases for any given special function. In the following we treat only those choices of $h$ such that equation (\ref{1}) allows us to derive new and interesting indefinite integrals. Last but not least, we checked that the software package Maple is unable to evaluate the indefinite integrals we computed in this paper.  

\subsection{The case of the confluent Heun function}\label{sch}
As in \cite{Dec} we restrict the local solutions to the analytic ones around the origin, here denoted by $H_c(\alpha,\beta,\gamma,\delta,\eta;x)$ but extensions are possible for other solutions around other singularities. The local solution is defined so that \cite{Dec}
\begin{equation}\label{incon}
H_c(\alpha,\beta,\gamma,\delta,\eta;0)=1,\quad
H^{'}_c(\alpha,\beta,\gamma,\delta,\eta;0)=\frac{(1+\gamma-\alpha)\beta+\gamma-\alpha+2\eta}{2(1+\beta)}.
\end{equation}
Some relations for the derivatives of confluent Heun functions have been derived in \cite{Ishk}. However, such formulae apply only to the  canonical form of the confluent Heun equation. Since Maple does not use this form, the validity of the aforementioned results cannot be checked using Maple and therefore, we need to derive corresponding formulae for when the confluent Heun equation is given as in Table~\ref{Heun}. It should also be mentioned that an additional formula for the derivatives of arbitrary order of a confluent Heun function has been obtained in \cite{Fiziev} under the condition that such a function reduces to a polynomial. This result will not be applied in the present work since we are mainly interested in deriving indefinite integrals for confluent Heun functions that are not necessarily represented by polynomials. The main idea in \cite{Ishk} is to differentiate once the confluent Heun equation and then, to subtract from it the confluent Heun equation multiplied by $q^{'}/q$. This leads to the following third order differential equation
\begin{equation}\label{dos}
u^{'''}(x)+\left[p(x)-\frac{q^{'}(x)}{q(x)}\right]u^{''}(x)+\left[p^{'}(x)+q(x)-\frac{q^{'}(x)}{q(x)}p(x)\right]u^{'}(x)=0
\end{equation}
with
\begin{equation}
p(x)-\frac{q^{'}(x)}{q(x)}=\alpha+\frac{2+\beta}{x}+\frac{2+\gamma}{x-1}-\frac{(\beta+\gamma+2)\alpha+2\delta}{[(\beta+\gamma+2)\alpha+2\delta]x-(\beta+1)\alpha+(\gamma+1)\beta+2\eta+\gamma}
\end{equation}
and
\[
p^{'}(x)+q(x)-\frac{q^{'}(x)}{q(x)}p(x)=-\frac{1+\beta}{x^2}-\frac{1+\gamma}{(x-1)^2}
\]
\[
+\frac{[(\beta+\gamma+2)\alpha+2\delta]x-(\beta+1)\alpha+(\gamma+1)\beta+2\eta+\gamma}{2x(x-1)}-\left(-\frac{1}{x}-\frac{1}{x-1}\right.
\]
\begin{equation}
\left.+\frac{(\beta+\gamma+2)\alpha+2\delta}{[(\beta+\gamma+2)\alpha+2\delta]x-(\beta+1)\alpha+(\gamma+1)\beta+2\eta+\gamma}\right)\left(\alpha+\frac{1+\beta}{x}+\frac{1+\gamma}{x-1}\right).
\end{equation}
Equation (\ref{dos}) represents a second order ODE for the first derivative of the Heun confluent function and it exhibits an additional singularity at the point
\begin{equation}
x_0=\frac{(\beta+1)\alpha-(\gamma+1)\beta-2\eta-\gamma}{(\beta+\gamma+2)\alpha+2\delta}.
\end{equation}
If the singularity at $x_0$ coincides with one of the singularities at $0$, $1$, $\infty$ the number of singularities will not change and in this case $0$ and $1$ remain regular singular points of (\ref{dos}) and the point at infinity continues to have the same rank. This will happen in the following three distinct cases
\begin{enumerate}
\item
$(\beta+\gamma+2)\alpha+2\delta=0$,
\item
$-(\beta+1)\alpha+(\gamma+1)\beta+2\eta+\gamma=0$,
\item
$(\beta+\gamma+2)\alpha+2\delta=(\beta+1)\alpha-(\gamma+1)\beta-2\eta-\gamma$ or equivalently
\begin{equation}
(\alpha+\beta)(\gamma+1)+2\delta=-2\eta-\gamma.
\end{equation}
\end{enumerate}
In each of the three cases above, equation (\ref{dos}) represents a confluent Heun equation for the first derivative of a confluent Heun function with modified parameters as compared to the confluent Heun equation given in Table~\ref{Heun}. Let us consider the first case above. It is straightforward to check with Maple that $u(x)=H_c\left(\alpha,\beta,\gamma,-(\beta+\gamma+2)\alpha/2,\eta;x\right)$ is a particular solution of (\ref{dos}) and 
\begin{equation}
w(x)=H_c\left(\alpha,\beta+1,\gamma+1,-(\beta+\gamma)\frac{\alpha}{2},\frac{\beta}{2}+\frac{\gamma}{2}-\frac{\alpha}{2}+\frac{1}{2}+\eta;x\right)
\end{equation}
a particular solution of
\[
w^{''}(x)+\left(\alpha+\frac{2+\beta}{x}+\frac{2+\gamma}{x-1}\right)w^{'}(x)+
\]
\begin{equation}\label{ww}
\frac{4\alpha x-(\beta+3)\alpha+(\gamma+3)\beta+2\eta+3\gamma+4}{2x(x-1)}w(x)=0,
\end{equation}
which is the counterpart of (\ref{dos}) when $(\beta+\gamma+2)\alpha+2\delta=0$. Furthermore, it can also be verified with Maple that $u^{'}(x)$ is a particular solution of (\ref{ww}). This implies that $u^{'}(x)=cw(x)$ with $c$ a proportionality constant that can be determined by means of the initial conditions (\ref{incon}). Taking into account that $w(0)=1$ we find
\[
H^{'}_c\left(\alpha,\beta,\gamma,-(\beta+\gamma+2)\frac{\alpha}{2},\eta;x\right)=\frac{(1+\gamma-\alpha)\beta+\gamma-\alpha+2\eta}{2(1+\beta)}\cdot
\]
\begin{equation}\label{der1}
\cdot H_c\left(\alpha,\beta+1,\gamma+1,-(\beta+\gamma)\frac{\alpha}{2},\frac{\beta}{2}+\frac{\gamma}{2}-\frac{\alpha}{2}+\frac{1}{2}+\eta;x\right).
\end{equation}
Concerning the second case  we choose
\begin{equation}
\eta=(\beta+1)\frac{\alpha}{2}-(\gamma+1)\frac{\beta}{2}-\frac{\gamma}{2}
\end{equation}
and (\ref{dos}) becomes
\begin{equation}
u^{'''}(x)+\left(\alpha+\frac{1+\beta}{x}+\frac{2+\gamma}{x-1}\right)u^{''}(x)+\frac{[(\beta+\gamma+4)\alpha+2\delta]x^2+2(\beta+1)}{2x^2(x-1)}u^{'}(x)=0.
\end{equation}
If we introduce the transformation $u^{'}(x)=x^s w(x)$ and choose the parameter $s$ so that the numerator going with the term $x^{-2}$ and appearing in the multiplicative coefficient of $w$ vanishes, then $s=-1$ or $s=-1-\beta$ and $w$ satisfies the confluent Heun equation
\[
w^{''}(x)+\left(\alpha+\frac{1+\beta+2s}{x}+\frac{2+\gamma}{x-1}\right)w^{'}(x)
\]
\begin{equation}
+\frac{[(\beta+\gamma+2s+4)\alpha+2\delta]x-2[(\alpha-\gamma-2)s-(\beta+1)]}{2x(x-1)}w(x)=0.
\end{equation}
Recalling that in the present case (\ref{incon}) gives $u^{'}(0)=0$ and proceeding as in the first case above yields
\[
H^{'}_c\left(\alpha,\beta,\gamma,\delta,(\beta+1)\frac{\alpha}{2}-(\gamma+1)\frac{\beta}{2}-\frac{\gamma}{2};x\right)=
\]
\begin{equation}\label{der2}
x^s H_c\left(\alpha,2s+\beta,\gamma+1,\frac{\alpha}{2}+\delta,(\alpha-\gamma)\frac{\beta}{2}+\frac{\alpha}{2}-\frac{\gamma}{2}+\frac{1}{2};x\right).
\end{equation}
Note that in the case $s=-1-\beta$ formula (\ref{der2}) requires that $\Re{(\beta)}<-1$ in order not to violate the initial condition $u^{'}(0)=0$. A formula similar to (\ref{der2}) can be obtained when $(\alpha+\beta)(\gamma+1)+2\delta=-2\eta-\gamma$. Taking into account that (\ref{f}) gives $f(x)=x^{1+\beta} (x-1)^{1+\gamma} e^{\alpha x}$,  if we choose $h$ according to (\ref{scelta_h}), we find
\[
\int\!x^{\beta+m-1}(x-1)^{\gamma}e^{\alpha x+\rho x^\ell}
\mathfrak{F}(x,k_1,k_2)
H_c(\alpha,\beta,\gamma,\delta,\eta;x)\,\mathrm{d}x= 
\]
\begin{equation}\label{Heun2}
2x^{\beta+m}(x-1)^{1+\gamma}e^{\alpha x+\rho x^\ell}
\mathfrak{G}(x,k_1,k_2)+c
\end{equation}
with 
\begin{eqnarray}
\mathfrak{F}(x,k_1,k_2)&=&\left\{
\begin{array}{c}
2x\mathfrak{p}_1(x,k_1)\cos{(k_1 x)}+\mathfrak{p}_2(x,k_1)\sin{(k_1 x)}\\
\mathfrak{p}_2(x,k_2)\cos{(k_2 x)}-2x\mathfrak{p}_1(x,k_2)\sin{(k_2 x)}
\end{array}
\right\},\\
\mathfrak{G}(x,k_1,k_2)&=&\left\{
\begin{array}{c}
\mathfrak{q}(x)\sin{(k_1 x)}+k_1 x\cos{(k_1 x)}H_c(\alpha,\beta,\gamma,\delta,\eta;x)\\
\mathfrak{q}(x)\cos{(k_2 x)}-k_2 x\sin{(k_2 x)}H_c(\alpha,\beta,\gamma,\delta,\eta;x)
\end{array}
\right\},\\
\mathfrak{q}(x)&=&(m+\rho\ell x^\ell)H_c(\alpha,\beta,\gamma,\delta,\eta;x)-xH^{'}_c(\alpha,\beta,\gamma,\delta,\eta;x),\\
\mathfrak{p}_1(x,k)&=&k\sum_{i=0}^2\mathfrak{a}_i x^i+2k\rho\ell x^\ell(x-1),\\
\mathfrak{p}_2(x,k)&=&\sum_{i=0}^3\mathfrak{b}_i x^i+2\rho\ell x^\ell\left[\sum_{i=0}^2\mathfrak{c}_i x^i+\rho\ell x^\ell(x-1)\right]
\end{eqnarray}
and
\begin{eqnarray}
\mathfrak{a_2}&=&\alpha,\quad \mathfrak{a}_1=\beta+\gamma-\alpha+2m+2,\quad \mathfrak{a}_0=-\beta-2m-1,\quad \mathfrak{b}_3=-2k^2,\\
\mathfrak{b}_2&=&(\beta+\gamma+2m+2)\alpha+2k^2+2\delta,\\
\mathfrak{b}_1&=&2m^2+2(\gamma-\alpha+\beta+1)m+(1+\beta)(\gamma-\alpha)+\beta+2\eta,\\
\mathfrak{b}_0&=&-2m(m+\beta),\quad
\mathfrak{c}_2=\alpha,\quad
\mathfrak{c}_1=\ell+\beta+\gamma-\alpha+2m+1,\\
\mathfrak{c}_0&=&-\ell-\beta-2m.
\end{eqnarray}
Note that the first derivative of the confluent Heun function appearing in the term $\mathfrak{q}(x)$ in (\ref{Heun2}) can be computed once the parameters of the Heun function have been chosen as in (\ref{der1}) or (\ref{der2}). Furthermore, by means of (\ref{zero}) we obtain immediately the following result
\begin{equation}\label{Hcnice}
\int\! x^{\beta}(x-1)^{\gamma}F(x)e^{\alpha x}H_c(\alpha,\beta,\gamma,\delta,\eta;x)\,\mathrm{d}x=
-2x^{1+\beta}(x-1)^{1+\gamma}e^{\alpha x}H^{'}_c(\alpha,\beta,\gamma,\delta,\eta;x)+c
\end{equation}
with
\begin{equation}
F(x)=\left[(\beta+\gamma+2)\alpha+2\delta\right]x-(\beta+1)\alpha+(\gamma+1)\beta+2\eta+\gamma.
\end{equation}
The derivative of the confluent Heun equation appearing in (\ref{Hcnice}) can be computed according to (\ref{der1}) or (\ref{der2}), if parameters are fixed accordingly. For instance, in the case corresponding to (\ref{der1}) we find
\[
\int\! x^{\beta}(x-1)^{\gamma}e^{\alpha x}H_c\left(\alpha,\beta,\gamma,-(\beta+\gamma+2)\frac{\alpha}{2},\eta;x\right)\,\mathrm{d}x=
\]
\begin{equation}
-\frac{1}{1+\beta}x^{1+\beta}(x-1)^{1+\gamma}e^{\alpha x}H_c\left(\alpha,\beta+1,\gamma+1,-(\beta+\gamma)\frac{\alpha}{2},\frac{\beta}{2}+\frac{\gamma}{2}-\frac{\alpha}{2}+\frac{1}{2}+\eta;x\right)+c.
\end{equation}
An integral involving $H_c(p,\gamma,\delta,\alpha,\sigma;x)$ alone can be constructed by making the choice $h(x)=1/\eta$ and $\alpha=\beta=\gamma=\delta=0$. In this case, $H_c(0,0,0,0,\eta;x)$ reduces to an hypergeometric function and the corresponding integral can be evaluated by Mathematica. Another possible choice is $h(x)=1/\delta$ while the parameters are taken to be $\alpha=\gamma=0$, $\beta=-1$, and $\eta=1/2$. Then, we obtain the integral
\begin{equation}\label{stanjel}
\int\! H_c(0,-1,0,\delta,1/2;x)\,\mathrm{d}x=\frac{1-x}{\delta}H^{'}_c(0,-1,0,\delta,1/2;x)+c.
\end{equation}
In this case the parameters satisfy the condition leading to formula (\ref{der2}) with $s=-1$. Hence, we are able to compute the derivative appearing in (\ref{stanjel}), and we get
\begin{equation}
\int\! H_c(0,-1,0,\delta,1/2;x)\,\mathrm{d}x=\frac{1-x}{\delta x}H_c(0,-3,1,\delta,1/2;x)+c.
\end{equation}
Moreover, we can also take $h$ to be a solution of the ODE (\ref{2}). In this case, an explicit analytic expression for the solution can be found only for certain choices of the parameters. For instance, for $\alpha=0$ we find $h(x)=x^{-\beta}{}_{2}F_{1}(-\beta,1+\gamma;1-\beta;x)$, and using (\ref{1}) together with 15.2.1 in \cite{abra} yields for $\beta\neq 1$
\[
\int\! (x-1)^{\gamma}U(x){}_{2}F_{1}(-\beta,1+\gamma;1-\beta;x)H_c(0,\beta,\gamma,\delta,\eta;x)\,\mathrm{d}x=
\]
\begin{equation}\label{nwfl}
2(x-1)^{1+\gamma}\left[\mathfrak{H}(x)H_c(0,\beta,\gamma,\delta,\eta;x)-x{}_{2}F_{1}(-\beta,1+\gamma;1-\beta;x)H^{'}_c(0,\beta,\gamma,\delta,\eta;x)\right]+c
\end{equation}
with
\begin{eqnarray}
U(x)&=&2\delta x+(\gamma+1)\beta+2\eta+\gamma,\\
\mathfrak{H}(x)&=&\beta\left[\frac{1+\gamma}{\beta-1}x{}_{2}F_{1}(1-\beta,2+\gamma;2-\beta;x)-{}_{2}F_{1}(-\beta,1+\gamma;1-\beta;x)\right].
\end{eqnarray}
Note that the derivative of the confluent Heun function in (\ref{nwfl}) can be computed by means of (\ref{der1}) with $\alpha=\delta=0$. We continue the analysis of the case of the confluent Heun function by observing that $h$ can also be taken to be a solution of the ODE (\ref{3}). In this case we find that
\begin{eqnarray}
h(x)&=&\mbox{exp}\left(-\int\!\frac{Mx+N}{Cx^2+2Bx+A}\,\mathrm{d}x\right),\label{K}\\
M&=&(\beta+\gamma+2)\alpha+2\delta,\quad N=-(\beta+1)\alpha+(\gamma+1)\beta+2\eta+\gamma,\\
C&=&2\alpha,\quad
B=\beta+\gamma+2-\alpha,\quad
A=-2(\beta+1).
\end{eqnarray}
Let $\Delta=AC-B^2$. Then, the  integral giving the function $h$ can be explicitly computed by means of $2.103.5$ in \cite{Grad} yielding
\begin{equation}
h(x)=\left\{
\begin{array}{ccc}
\left[Cx^2+2Bx+A\right]^{-\frac{M}{2C}}\mbox{exp}\left[\frac{MB-NC}{C\sqrt{\Delta}}\arctan{\left(\frac{Cx+B}{\sqrt{\Delta}}\right)}\right]&\mbox{if}~\Delta>0,\\
(Cx+B)^{-\frac{M}{C}}\mbox{exp}\left[\frac{NC-MB}{C(Cx+B)}\right]  & \mbox{if}~\Delta=0,\\
\left[Cx^2+2Bx+A\right]^{-\frac{M}{2C}}\left[\frac{Cx+B-\sqrt{-\Delta}}{Cx+B+\sqrt{-\Delta}}\right]^\frac{MB-NC}{2C\sqrt{-\Delta}}&\mbox{if}~\Delta<0.
\end{array}
\right.
\end{equation}
Since there is a singularity at $x_0=(\alpha-\beta-\gamma-2)/2\alpha$ in the case $\Delta=0$, we need to introduce some additional constraints on the parameters ensuring that $x_0$ lies outside the interval where the local solution of the confluent Heun equation is defined. Finally, (\ref{1}) leads to the following indefinite integral
\[
\int\!x^{1+\beta}(x-1)^{1+\gamma}e^{\alpha x}u(x)h(x)H_c(\alpha,\beta,\gamma,\delta,\eta;x)\,\mathrm{d}x=
\]
\begin{equation}\label{ciccio}
-x^{1+\beta}(x-1)^{1+\gamma} e^{\alpha x} h(x)\left[\frac{Mx+N}{Cx^2+2Bx+A}H_c(\alpha,\beta,\gamma,\delta,\eta;x)+H^{'}_c(\alpha,\beta,\gamma,\delta,\eta;x)\right]+c
\end{equation}
with
\begin{equation}
u(x)=\frac{M(C+M)x^2+2N(C+M)x+N^2+2BN-AM}{(Cx^2+2Bx+A)^2}.
\end{equation}
The derivative of the confluent Heun equation appearing in (\ref{ciccio}) can be computed according to (\ref{der1}) or (\ref{der2}), if the parameters are fixed accordingly. Further indefinite integrals can be obtained by choosing  the function $h$ to be a particular solution to (\ref{4}). An interesting case emerges when we choose the parameter $\eta$ so that it takes the value
\begin{equation}
\eta_0=(1+\beta)\frac{\alpha}{2}-(1+\gamma)\frac{\beta}{2}-\frac{\gamma}{2}.
\end{equation}
Then, we can express $h$ in terms of Bessel functions of the first and second kind as follows
\begin{equation}
h(x)=\sqrt{x-1}\left\{
\begin{array}{cc}
J_1(\Omega\sqrt{x-1})\\
Y_1(\Omega\sqrt{x-1})
\end{array}
\right\},\quad \Omega=\sqrt{2\alpha(\beta+\gamma+2)+4\delta}.
\end{equation}
Moreover, it can be checked that
\begin{equation}
h^{'}(x)=\frac{\Omega}{2}\left\{
\begin{array}{cc}
J_0(\Omega\sqrt{x-1})\\
Y_0(\Omega\sqrt{x-1})
\end{array}
\right\},
\end{equation}
where we used 9.1.27 in \cite{abra}. Then, by means of (\ref{1}) we obtain
\[
\int\!x^{\beta}(x-1)^{\gamma}e^{\alpha x}K(x)\left\{
\begin{array}{cc}
J_0(z(x))\\
Y_0(z(x))
\end{array}
\right\}H_c(\alpha,\beta,\gamma,\delta,\eta_0;x)\,\mathrm{d}x=x^{1+\beta}(x-1)^{1+\gamma}e^{\alpha x}\cdot
\]
\begin{equation}\label{hahaz}
\left[\left\{
\begin{array}{cc}
J_0(z(x))\\
Y_0(z(x))
\end{array}
\right\}H_c(\alpha,\beta,\gamma,\delta,\eta_0;x)-\frac{2}{\Omega}\sqrt{x-1}\left\{
\begin{array}{cc}
J_1(z(x))\\
Y_1(z(x))
\end{array}
\right\}H^{'}_c(\alpha,\beta,\gamma,\delta,\eta_0;x)\right]+c
\end{equation}
with $z(x)=\Omega\sqrt{x-1}$, $\Omega=\sqrt{2\alpha(\beta+\gamma+2)+4\delta}$, and
\begin{equation}
K(x)=\alpha x^2+(2+\beta+\gamma-\alpha)x-\beta-1.
\end{equation}
The derivative of the confluent Heun function in (\ref{hahaz}) can be evaluated with the help of (\ref{der1}) or (\ref{der2}) after the parameters are chosen accordingly. Furthermore, let us consider the particular solution $y(x)=H_c(\alpha,\beta,\gamma,\delta,\eta;x)$ to the confluent Heun equation and a conjugate ODE to the same equation with $\overline{q}(x)$ given by $q(x)$ as in Table~\ref{Heun} but with the parameter $\eta$ replaced by $-\eta$. For the conjugate ODE we pick the particular solution $h(x)=H_c(\alpha,\beta,\gamma,\delta,-\eta;x)$. Then, (\ref{intcon}) yields the following indefinite integral involving products of Heun functions
\[
\int\!x^{\beta}(x-1)^{\gamma}e^{\alpha x}H_c(\alpha,\beta,\gamma,\delta,\eta;x)H_c(\alpha,\beta,\gamma,\delta,-\eta;x)\,\mathrm{d}x=\frac{x^{1+\beta}}{2\eta}(x-1)^{1+\gamma}e^{\alpha x}\cdot
\]
\begin{equation}\label{hhh1}
\left[H^{'}_c(\alpha,\beta,\gamma,\delta,-\eta;x)H_c(\alpha,\beta,\gamma,\delta,\eta;x)-H_c(\alpha,\beta,\gamma,\delta,-\eta;x)H^{'}_c(\alpha,\beta,\gamma,\delta,\eta;x)\right]+c.
\end{equation}
provided that $\eta\neq 0$. The derivatives of the Heun functions in (\ref{hhh1}) can be computed according to (\ref{der1}) or (\ref{der2}) for an appropriate choice of the parameters.

\subsection{The case of the biconfluent Heun function}
As in \cite{us} we restrict the local solutions to the analytic ones around the origin, here denoted by $H_b(\alpha,\beta,\gamma,\delta;x)$. Furthermore, the local solution is defined so that \cite{HR}
\begin{equation}\label{anfang}
H_b(\alpha,\beta,\gamma,\delta;0)=1,\quad
H_b^{'}(\alpha,0,\gamma,0;0)=\frac{\delta+\beta(\alpha+1)}{2(1+\alpha)}.
\end{equation}
In the following, we will make use of some useful relations for the derivative of confluent Heun functions derived in \cite{HR} and summarized here below
\begin{eqnarray}
H_b(\alpha,0,\gamma,0;x)&=&{}_{1}F_{1}\left(\frac{\alpha+2-\gamma}{4},1+\frac{\alpha}{2};x^2\right),\label{ST1}\\
H_b^{'}(\alpha,0,\gamma,0;x)&=&\frac{(\alpha+2-\gamma)x}{\alpha+2}{}_{1}F_{1}\left(\frac{\alpha+6-\gamma}{4},2+\frac{\alpha}{2};x^2\right),\label{Ab1}
\end{eqnarray}
where ${}_{1}F_{1}(\cdot)$ denotes the confluent hypergeometric function. Using the same method presented at the beginning of Section~\ref{sch} together with (\ref{anfang}) it is possible to derive a new formula for the derivative of the biconfluent Heun function, namely
\begin{equation}\label{Ab2}
H^{'}_b(\alpha,\beta,\alpha+2,\delta;x)=\frac{\delta+\beta(\alpha+1)}{2(1+\alpha)}H_b(\alpha+1,\beta,\alpha-1,\beta+\delta;x),\quad\alpha\neq-1.
\end{equation}
Moreover, in the present case (\ref{f}) gives $f(x)=x^{\alpha+1}e^{-x^2-\beta x}$. If we choose $h$ according to (\ref{scelta_h}), we find by means of (\ref{1})
\[
\int\!x^{\alpha+m-1}e^{-x^2-\beta x+\rho x^\ell}
\widehat{\mathfrak{F}}(x,k_1,k_2)
H_b(\alpha,\beta,\gamma,\delta;x)\,\mathrm{d}x=
\]
\begin{equation}\label{HeunBC2}
x^{\alpha+m}e^{-x^2-\beta x+\rho x^\ell}
\widehat{\mathfrak{G}}(x,k_1,k_2)+c
\end{equation}
with 
\begin{eqnarray}
\widehat{\mathfrak{F}}(x,k_1,k_2)&=&\left\{
\begin{array}{c}
x\widehat{\mathfrak{p}}_1(x,k_1)\cos{(k_1 x)}+\widehat{\mathfrak{p}}_2(x,k_1)\sin{(k_1 x)}\\
\widehat{\mathfrak{p}}_2(x,k_2)\cos{(k_2 x)}-x\widehat{\mathfrak{p}}_1(x,k_2)\sin{(k_2 x)}
\end{array}
\right\},\\
\widehat{\mathfrak{G}}(x,k_1,k_2)&=&\left\{
\begin{array}{c}
\widehat{\mathfrak{q}}(x)\sin{(k_1 x)}+k_1 x\cos{(k_1 x)}H_b(\alpha,\beta,\gamma,\delta;x)\\
\widehat{\mathfrak{q}}(x)\cos{(k_2 x)}-k_2 x\sin{(k_2 x)}H_b(\alpha,\beta,\gamma,\delta;x)
\end{array}
\right\},\\
\widehat{\mathfrak{q}}(x)&=&(m+\rho\ell x^\ell)H_b(\alpha,\beta,\gamma,\delta;x)-xH^{'}_b(\alpha,\beta,\gamma,\delta;x),\label{qk}\\
\widehat{\mathfrak{p}}_1(x,k)&=&k\sum_{i=0}^2\widehat{\mathfrak{a}}_i x^i+2k\rho\ell x^\ell,\\
\widehat{\mathfrak{p}}_2(x,k)&=&\sum_{i=0}^2\widehat{\mathfrak{b}}_i x^i+\rho\ell x^\ell\left[\sum_{i=0}^2\widehat{\mathfrak{c}}_i x^i+\rho\ell x^\ell\right]
\end{eqnarray}
and
\begin{eqnarray}
\widehat{\mathfrak{a}}_2&=&\widehat{\mathfrak{c}}_2=-2,~ \widehat{\mathfrak{a}}_1=\widehat{\mathfrak{c}}_1=-\beta,~ \widehat{\mathfrak{a}}_0=1+2m+\alpha\\
\widehat{\mathfrak{b}}_2&=&-k^2-\alpha+\gamma-2m-2,~\widehat{\mathfrak{b}}_1=-\frac{\beta(\alpha+2m+1)+\delta}{2},\\
\widehat{\mathfrak{b}}_0&=&m(m+\alpha),~
\widehat{\mathfrak{c}_0}=\ell+\alpha+2m.
\end{eqnarray}
Note that for an appropriate choice of the parameters of the biconfluent Heun function the derivative in (\ref{HeunBC2}) can be computed by means of (\ref{Ab1}), or (\ref{Ab2}). Furthermore, by means of (\ref{zero}) we obtain immediately the following result
\begin{equation}\label{HDCnice}
\int\! x^{\alpha}(A_1 x-A_2)e^{-x^2-\beta x}H_b(\alpha,\beta,\gamma,\delta;x)\,\mathrm{d}x=-2x^{\alpha+1}e^{-x^2-\beta x}H_b^{'}(\alpha,\beta,\gamma,\delta;x)+c.
\end{equation}
with $A_1=2(\gamma-\alpha-2)$ and $A_2=\delta+\beta(\alpha+1)$. If $\gamma=\alpha+2$, the derivative of the biconfluent Heun function in (\ref{HDCnice}) can be computed with the help of (\ref{Ab2}), and in this case we find
\begin{equation}\label{HDCspc}
\int\! x^{\alpha}e^{-x^2-\beta x}H_b(\alpha,\beta,\alpha+2,\delta;x)\,\mathrm{d}x=\frac{x^{\alpha+1}}{\alpha+1}e^{-x^2-\beta x}H_b(\alpha+1,\beta,\alpha-1,\beta+\delta;x)+c
\end{equation}
provided that $\alpha\neq 1$. Moreover, we can also take $h$ to be a solution of the ODE (\ref{2}). In the most general case, a particular solution can be formally written as
\begin{equation}
h(x)=\int\! x^{-\alpha-1}e^{x^2+\beta x}\,\mathrm{d}x.
\end{equation}
However, the above integral can be explicitly solved only for some special choices of the parameters. For instance, if $\alpha=-1$, 2.325(13) in \cite{Grad} gives
\begin{equation}
h(x)=\frac{\sqrt{\pi}}{2}e^{-\frac{\beta^2}{4}}\mbox{erfi}\left(x+\frac{\beta}{2}\right),
\end{equation}
where $\mbox{erfi}(\cdot)$ is the imaginary error function defined as $\mbox{erfi}(z)=-i\mbox{erf}(iz)$ with $\mbox{erf}$ denoting the error function. Using 7.1.1 in \cite{abra} and (\ref{1}) yield 
\[
\int\! e^{-x^2-\beta x}\left(\gamma-1-\frac{\delta}{2x}\right)\mbox{erfi}\left(x+\frac{\beta}{2}\right)H_b(-1,\beta,\gamma,\delta;x)\,\mathrm{d}x=
\]
\begin{equation}\label{HDChp}
\frac{2e^{\frac{\beta^2}{4}}}{\sqrt{\pi}}H_b(-1,\beta,\gamma,\delta;x)-e^{-x^2-\beta x}\mbox{erfi}\left(x+\frac{\beta}{2}\right)H_b^{'}(-1,\beta,\gamma,\delta;x)+c.
\end{equation}
Note that in the case $\alpha=-1$ only (\ref{Ab1}) can be used to evaluate the derivative of the biconfluent Heun function. We continue the analysis of the case of the biconfluent Heun function by observing that $h$ can also be taken to be a solution of the ODE (\ref{3}). If $\Delta=-2(\alpha+1)-(\beta^2/4)$, the  integral giving the function $h$ can be explicitly computed by means of $2.103.5$ in \cite{Grad} yielding
\begin{equation}
h(x)=\left\{
\begin{array}{ccc}
(2x^2+\beta x-\alpha-1)^{\frac{\gamma-\alpha-2}{4}}\mbox{exp}\left[-\frac{2\delta+\beta(\alpha+\gamma)}{4\sqrt{\Delta}}\arctan{\left(\frac{4x+\beta}{2\sqrt{\Delta}}\right)}\right]&\mbox{if}~\Delta>0,\\
\left(x+\frac{\beta}{4}\right)^{\frac{\gamma-\alpha-2}{2}}\mbox{exp}\left(\frac{\beta(\gamma+\alpha)+2\delta}{2(4x+\beta)}\right)  & \mbox{if}~\Delta=0,\\
(2x^2+\beta x-\alpha-1)^{\frac{\gamma-\alpha-2}{4}}\left[\frac{4x+\beta+2\sqrt{-\Delta}}{4x+\beta-2\sqrt{-\Delta}}\right]^\frac{2\delta+\beta(\alpha+\gamma)}{8\sqrt{-\Delta}}&\mbox{if}~\Delta<0.
\end{array}
\right.
\end{equation}
At this point a couple of remarks are in order. First of all, observe that the condition $\Delta>0$ requires that $\alpha<-1$. Moreover, if $\Delta=0$, we have $\alpha=-1-(\beta^2/8)$, and the corresponding solution has a singularity at $x_0=-\beta/4$. In this case, we need to require that that $x_0$ lies outside the interval where the local solution of the biconfluent Heun equation is defined. This observation leads to the constraint $\beta\in\mathbb{R}\backslash[-4,0]$. Finally, (\ref{1}) leads to the following indefinite integral
\[
\int\!x^{\alpha+1} e^{-x^2-\beta x}U(x)h(x)H_b(\alpha,\beta,\gamma,\delta;x)\,\mathrm{d}x=
\]
\begin{equation}\label{ciccio1}
2x^{\alpha+1} e^{-x^2-\beta x} h(x)\left[S(x)H_b(\alpha,\beta,\gamma,\delta;x)-2H^{'}_b(\alpha,\beta,\gamma,\delta;x)\right]+c
\end{equation}
with
\begin{equation}\label{Mx}
U(x)=\frac{A_2 x^2+A_1 x+A_0}{(2x^2+\beta x-\alpha-1)^2},\quad S(x)=\frac{2(\gamma-\alpha-2)x-\beta(\alpha+1)-\delta}{2x^2+\beta x-\alpha-1},
\end{equation}
and coefficients
\begin{eqnarray}
A_2&=&4\left[(\alpha-\gamma)^2+2(3\alpha-3\gamma+4)\right],\\
A_1&=&4\left[\alpha\beta(\alpha-\gamma+5)+\alpha\delta+(\beta+\delta)(4-\gamma)\right],\\
A_0&=&\alpha\beta(\alpha\beta+4\beta+2\delta)+\delta^2+3\beta^2+4\left[\alpha^2+\beta\delta-\alpha\gamma+3\alpha-\gamma+2\right].
\end{eqnarray}
The derivative of the confluent Heun equation appearing in (\ref{ciccio1}) can be computed according to (\ref{Ab1}), and (\ref{Ab2}) if the parameters are fixed accordingly. Furthermore, we can consider the particular solution $y(x)=H_b(\alpha,\beta,\gamma,\delta;x)$ to the biconfluent Heun equation and an ODE conjugate to the same equation with $\overline{q}(x)$ given by the corresponding $q(x)$ in Table~\ref{Heun} with $\gamma$ and $\delta$ replaced by $-\gamma$ and $-\delta$, respectively. Then, (\ref{intcon}) yields the following indefinite integral involving products of Heun functions
\[
\int\!x^{\alpha}(2\gamma x-\delta)e^{-x^2-\beta x}H_b(\alpha,\beta,\gamma,\delta;x)H_b(\alpha,\beta,-\gamma,-\delta;x)\,\mathrm{d}x=x^{\alpha+1}e^{-x^2-\beta x}\cdot
\]
\begin{equation}\label{hhhX}
\left[H_b(\alpha,\beta,\gamma,\delta;x)H_b^{'}(\alpha,\beta,-\gamma,-\delta;x)-H_b^{'}(\alpha,\beta,\gamma,\delta;x)H_b(\alpha,\beta,-\gamma,-\delta;x)\right]+c.
\end{equation}

\subsection{The case of the doubly confluent Heun function} 
Let $H_d(\alpha,\beta,\gamma,\delta;x)$ denote the restriction of the local solutions to the analytic ones around the origin satisfying the initial conditions $H_d(\alpha,\beta,\gamma,\delta;0)=1$ and $H^{'}_d(\alpha,\beta,\gamma,\delta;0)=0$ \cite{Dec}. Moreover, in the present case (\ref{f}) gives $f(x)=(x^2-1)e^{\frac{\alpha x}{x^2-1}}$. To the best of our knowledge, no formulae for the first derivative of the doubly confluent Heun function seem to be available in the literature. If we choose $h$ according to (\ref{scelta_h}), we get by means of (\ref{1})
\begin{equation}\label{HeunDC2}
\int\!\frac{x^{m-2}e^{\frac{\alpha x}{x^2-1}+\rho x^\ell}}{(x^2-1)^2}
\widetilde{\mathfrak{F}}(x,k_1,k_2)
H_d(\alpha,\beta,\gamma,\delta;x)\,\mathrm{d}x=
x^{m-1}(x^2-1)e^{\frac{\alpha x}{x^2-1}+\rho x^\ell}
\widetilde{\mathfrak{G}}(x,k_1,k_2)+c
\end{equation}
with 
\begin{eqnarray}
\widetilde{\mathfrak{F}}(x,k_1,k_2)&=&\left\{
\begin{array}{c}
x\widetilde{\mathfrak{p}}_1(x,k_1)\cos{(k_1 x)}+\widetilde{\mathfrak{p}}_2(x,k_1)\sin{(k_1 x)}\\
\widetilde{\mathfrak{p}}_2(x,k_2)\cos{(k_2 x)}-x\widetilde{\mathfrak{p}}_1(x,k_2)\sin{(k_2 x)}
\end{array}
\right\},\\
\widetilde{\mathfrak{G}}(x,k_1,k_2)&=&\left\{
\begin{array}{c}
\widetilde{\mathfrak{q}}(x)\sin{(k_1 x)}+k_1 x\cos{(k_1 x)}H_d(\alpha,\beta,\gamma,\delta;x)\\
\widetilde{\mathfrak{q}}(x)\cos{(k_2 x)}-k_2 x\sin{(k_2 x)}H_d(\alpha,\beta,\gamma,\delta;x)
\end{array}
\right\},\\
\widetilde{\mathfrak{q}}(x)&=&(m+\rho\ell x^\ell)H_d(\alpha,\beta,\gamma,\delta;x)-xH^{'}_d(\alpha,\beta,\gamma,\delta;x),\label{qk1}\\
\widetilde{\mathfrak{p}}_1(x,k)&=&k\sum_{i=0}^6\widetilde{\mathfrak{a}}_i x^i+2k\rho\ell x^{\ell}(x^2-1)^3,\\
\widetilde{\mathfrak{p}}_2(x,k)&=&\sum_{i=0}^8\widetilde{\mathfrak{b}}_i x^i+\rho\ell x^\ell\left[\sum_{i=0}^6\widetilde{\mathfrak{c}}_i x^i+\rho\ell x^{\ell}(x^2-1)^3\right],
\end{eqnarray}
and
\begin{eqnarray}
\widetilde{\mathfrak{a}}_6&=&2m+2,~\widetilde{\mathfrak{a}}_5=-\widetilde{\mathfrak{a}}_1=-\alpha,~\widetilde{\mathfrak{a}}_4=-6m-4,~\widetilde{\mathfrak{a}}_3=0,~
\widetilde{\mathfrak{a}}_2=6m+2,\\
\widetilde{\mathfrak{a}}_0&=&-2m,~\widetilde{\mathfrak{b}}_8=-k^2,~\widetilde{\mathfrak{b}}_7=0,~\widetilde{\mathfrak{b}}_6=3k^2+m^2+m,~\widetilde{\mathfrak{b}}_5=-\widetilde{\mathfrak{b}}_1=-\alpha m,\\
\widetilde{\mathfrak{b}}_4&=&\beta-\widetilde{\mathfrak{b}}_6,~\widetilde{\mathfrak{b}}_3=\gamma+2\alpha,~\widetilde{\mathfrak{b}}_2=k^2+3m^2-m+\delta,~\widetilde{\mathfrak{b}}_0=m-m^2,\\
\widetilde{\mathfrak{c}}_6&=&\ell+2m+1,~\widetilde{\mathfrak{c}}_5=-\widetilde{\mathfrak{c}}_1=-\alpha,~\widetilde{\mathfrak{c}}_4=2-3\widetilde{\mathfrak{c}}_6,~\widetilde{\mathfrak{c}}_3=0,\\
\widetilde{\mathfrak{c}}_2&=&3\widetilde{\mathfrak{c}}_6-4,~\widetilde{\mathfrak{c}}_0=2-\widetilde{\mathfrak{c}}_6.
\end{eqnarray}
Furthermore, by means of (\ref{zero}) we obtain immediately the following result
\begin{equation}\label{dc1}
\int\!\frac{e^{\frac{\alpha x}{x^2-1}}}{(x^2-1)^2}S(x)H_d(\alpha,\beta,\gamma,\delta;x)\,\mathrm{d}x=
-(x^2-1)e^{\frac{\alpha x}{x^2-1}}H^{'}_d(\alpha,\beta,\gamma,\delta;x)+c
\end{equation}
with 
\begin{equation}
S(x)=\beta x^2+(\gamma+2\alpha)x+\delta.
\end{equation}
Moreover, if we choose the function $h$ to be a solution of (\ref{2}) with $\alpha=0$, namely 
\begin{equation}
h(x)=\ln{\sqrt{\frac{x-1}{x+1}}},
\end{equation}
equation (\ref{1}) gives
\[
\int\!\frac{\beta x^2+\gamma x+\delta}{(x^2-1)^2}\ln{\sqrt{\frac{x-1}{x+1}}}H_d(0,\beta,\gamma,\delta;x)\,\mathrm{d}x=
\]
\begin{equation}
H_d(0,\beta,\gamma,\delta;x)-(x^2-1)\ln{\sqrt{\frac{x-1}{x+1}}}H^{'}_d(0,\beta,\gamma,\delta;x).
\end{equation}
We can also derive a new indefinite integral by choosing  the function $h$ to be a particular solution to (\ref{3}). In the special case $\alpha=0$ we have
\begin{equation}
h(x)=x^{-\frac{\delta}{2}}(x-1)^\frac{2\delta+\gamma}{8}(x+1)^\frac{2\delta-\gamma}{8}e^\frac{\gamma x+\beta+\delta}{4(x^2-1)},
\end{equation}
and (\ref{1}) gives
\[
\int\! x^{-2-\frac{\delta}{2}}(x-1)^{\frac{\delta}{4}+\frac{\gamma}{8}-3}(x+1)^{\frac{\delta}{4}-\frac{\gamma}{8}-3}e^\frac{\gamma x+\beta+\delta}{4(x^2-1)}\mathfrak{U}(x)H_d(0,\beta,\gamma,\delta;x)\,\mathrm{d}x=
\]
\begin{equation}
-h(x)\left[\frac{2(\beta x^2+\gamma x+\delta)}{x(x^2-1)}H_d(0,\beta,\gamma,\delta;x)+4(x^2-1)H^{'}_d(0,\beta,\gamma,\delta;x)\right]
\end{equation}
with
\[
\mathfrak{U}(x)=6\beta x^6+8\gamma x^5+(\beta^2-4\beta+10\delta)x^4+2\gamma(\beta-4)x^3+
\]
\begin{equation}
(2\beta\delta+\gamma^2-2\beta-12\delta)x^2+2\gamma\delta x+\delta^2+2\delta.
\end{equation}
Finally, we consider the particular solution $y(x)=H_d(\alpha,\beta,\gamma,\delta;x)$ to the doubly confluent Heun equation and an ODE conjugate to the same equation with $\overline{q}(x)$ given for example by the corresponding $q(x)$ in Table~\ref{Heun} with $\gamma$ replaced by $-\gamma$. Then, (\ref{intcon}) yields the following indefinite integral for $\gamma\neq 0$
\[
\int\!\frac{x}{(x^2-1)^2}e^{\frac{\alpha x}{x^2-1}}H_d(\alpha,\beta,\gamma,\delta;x)H_d(\alpha,\beta,-\gamma,\delta;x)\,\mathrm{d}x= 
\]
\begin{equation}\label{hhhXc}
\frac{x^2-1}{2\gamma}e^{\frac{\alpha x}{x^2-1}}W\left(H_d(\alpha,\beta,\gamma,\delta;x),H_d(\alpha,\beta,-\gamma,\delta;x);x)\right)+c,
\end{equation}
where $W$ denotes the Wronskian.

\subsection{The case of the triconfluent Heun function}
We restrict the local solutions to the analytic one around the origin, here denoted by $H_t(\alpha,\beta,\gamma;x)$, and satisfying the initial conditions $H_t(\alpha,\beta,\gamma;0)=1$ and $H^{'}_t(\alpha,\beta,\gamma;0)=0$ \cite{Dec}. In the present case,  (\ref{f}) gives $f(x)=e^{-x^3-\gamma x}$. If we choose $h$ according to (\ref{scelta_h}), we find by means of (\ref{1})
\[
\int\!x^{m-2}e^{-x^3-\gamma x+\rho x^\ell}
\left\{
\begin{array}{c}
x\mathfrak{P}_1(x,k_1)\cos{(k_1 x)}+\mathfrak{P}_2(x,k_1)\sin{(k_1 x)}\\
\mathfrak{P}_2(x,k_2)\cos{(k_2 x)}-x\mathfrak{P}_1(x,k_2)\sin{(k_2 x)}
\end{array}
\right\}
H_t(\alpha,\beta,\gamma;x)\,\mathrm{d}x=
\]
\begin{equation}\label{Heuntri}
x^{m-1}e^{-x^3-\gamma x+\rho x^\ell}
\left\{
\begin{array}{c}
\mathfrak{Q}(x)\sin{(k_1 x)}+k_1 x\cos{(k_1 x)}H_t(\alpha,\beta,\gamma;x)\\
\mathfrak{Q}(x)\cos{(k_2 x)}-k_2 x\sin{(k_2 x)}H_t(\alpha,\beta,\gamma;x)
\end{array}
\right\}+c
\end{equation}
with 
\begin{eqnarray}
\mathfrak{Q}(x)&=&(m+\rho\ell x^\ell)H_t(\alpha,\beta,\gamma;x)-xH^{'}_t(\alpha,\beta,\gamma;x),\\
\mathfrak{P}_1(x,k)&=&-k(3x^3+\gamma x-2m)+2k\rho\ell x^\ell,\\
\mathfrak{P}_2(x,k)&=&\sum_{i=0}^3\mathfrak{s}_i x^i+\rho\ell x^\ell\left[\sum_{i=0}^3\mathfrak{t}_i x^i+\rho\ell x^\ell\right],
\end{eqnarray}
and
\begin{eqnarray}
\mathfrak{s}_3&=&\beta-3m-3,\quad\mathfrak{s}_2=\alpha-k^2,\quad\mathfrak{s}_1=-\gamma m,\quad\mathfrak{s}_0=m(m-1),\\
\mathfrak{t}_3&=&-3,\quad
\mathfrak{t}_2=0,\quad\mathfrak{t}_1=-\gamma,\quad\mathfrak{t}_0=\ell+2m-1.
\end{eqnarray}
We can also take $h$ to be a solution of (\ref{2}) with $\gamma=0$. In this case, we find that a particular solution can be expressed in terms of the upper incomplete Gamma function as follows
\begin{equation}\label{h1}
h(x)=\Gamma\left(\frac{1}{3},-x^3\right),
\end{equation}
where we used 2.325(6) in \cite{Grad}. At this point (\ref{1}) yields
\[
\int\left[\alpha+(\beta-3)x\right]e^{-x^3}\Gamma\left(\frac{1}{3},-x^3\right)H_t(\alpha,\beta,0;x)\,\mathrm{d}x=
\]
\begin{equation}
\frac{3x^2}{(-x^3)^\frac{2}{3}}H_t(\alpha,\beta,0;x)-e^{-x^3}\Gamma\left(\frac{1}{3},-x^3\right)H^{'}_t(\alpha,\beta,0;x),
\end{equation}
where the derivative of the upper incomplete Gamma function has been computed with 6.5.25 in\cite{abra}. Furthermore, it is also possible to take $h$ to be a solution of (\ref{3}). In this case, 2.103.5 in \cite{Grad} gives
\begin{equation}\label{aidusina}
h(x)=\left\{ 
\begin{array}{lll}
         (3x^2+\gamma)^{\frac{\beta-3}{6}}e^{\frac{\alpha}{\sqrt{3\gamma}}\arctan{\left(\frac{\sqrt{3}x}{\sqrt{\gamma}}\right)}} & \mbox{if}~\gamma>0, \\
         (3x^2+\gamma)^{\frac{\beta-3}{6}}\left(\frac{3x-\sqrt{-3\gamma}}{3x+\sqrt{-3\gamma}}\right)^{\frac{\alpha}{2\sqrt{-3\gamma}}} & \mbox{if}~\gamma<0,\\
		x^{\frac{\beta-3}{3}}e^{-\frac{\alpha}{3x}}  & \mbox{if}~\gamma=0. 
\end{array} \right.
\end{equation}
At this point we can use (\ref{1}) to get
\[
\int e^{-x^3-\gamma x}\frac{\mathfrak{K}(x)h(x)}{(3x^2+\gamma)^2}H_t(\alpha,\beta,\gamma;x)\,\mathrm{d}x=
\]
\begin{equation}
e^{-x^3-\gamma x}h(x)\left[\frac{(\beta-3)x+\alpha}{3x^2+\gamma}T_\ell(\alpha,\beta,\gamma;x)-T^{'}_\ell(\alpha,\beta,\gamma;x)\right]
\end{equation}
with $h$ given as in (\ref{aidusina}) and
\begin{equation}
\mathfrak{K}(x)=(\beta-3)(\beta-6)x^2+2\alpha(\beta-6)x+\alpha^2+\gamma(\beta-3).
\end{equation}
As a last case, let us take a particular solution $y(x)=H_t(\alpha,\beta,\gamma;x)$ of the triconfluent Heun function and a conjugate ODE with $\overline{q}(x)$ obtained from $q(x)$ under the transformation $\alpha\to-\alpha$ while all other parameters remain unchanged. Then, if $\alpha\neq 0, $we immediately obtain from (\ref{intcon})
\[
\int e^{-x^3-\gamma x}H_t(-\alpha,\beta,\gamma;x)H_t(\alpha,\beta,\gamma;x)\,\mathrm{d}x=
\]
\begin{equation}
\frac{e^{-x^3-\gamma x}}{2\alpha}\left[H^{'}_t(-\alpha,\beta,\gamma;x)H_t(\alpha,\beta,\gamma;x)-H_t(-\alpha,\beta,\gamma;x)H^{'}_t(\alpha,\beta,\gamma;x)\right]+c.
\end{equation}
A similar result can be obtained by considering a conjugate ODE constructed from $q(x)$ keeping all parameters unchanged except $\beta\to-\beta$.

\section{Comments and conclusions}
We applied the so-called Lagrangian method to obtain indefinite integrals of functions belonging to the family of Heun confluent functions. This approach allowed us to derive several novel indefinite integrals for the confluent, biconfluent, doubly confluent, and triconfluent Heun functions for which sample results have been provided. Our findings only scratch the surface of the wealth of new integral formulae one may obtain by using the aforementioned method. Some of the results presented still involve derivatives of the confluent Heun functions. This is mainly due to the fact that the available formulae for derivatives of confluent Heun functions in the literature are scarce and moreover, they hold for very special choices of the parameters entering in the corresponding ODE. Some preliminary results allowing us to compute the first derivative of confluent and biconfluent Heun functions have been obtained under the minimal assumption of fixing only one parameter. More general results in this direction will be presented separately.

\end{document}